\documentclass[12pt]{amsart}
\usepackage{tikz,
  fullpage,
  graphicx,amssymb,epic,eepic,cancel,
color,hyphenat,
  multicol, 
}
\usepackage[all]{xy}

\theoremstyle{definition}

\def\R{{\mathbb R}}

\def\Q{{\mathbb Q}}

\def\calK{{\mathcal K}}
\def\calB{{\mathcal B}}

\def\-{\setminus}

\begin{document}


\title{What is a Singular Knot?}

\author{Zsuzsanna Dancso}
\thanks{}
\address{Zsuzsanna Dancso,
  School of Mathematics and Statistics\\
  The Universit of Sydney\\
  Sydney NSW 2006\\
  Australia
}
\email{zsuzsanna.dancso@sydney.edu.au}

\begin{abstract} A singular knot is an immersed circle in $\R^{3}$ with finitely many transverse double points. The study of singular knots was initially motivated by the study of Vassiliev invariants. Namely, singular knots give rise to a decreasing filtration on the infinite dimensional vector space spanned by isotopy types of knots: this is called the Vassiliev filtration, and the study of the corresponding associated graded space has lead to many insights in knot theory. The Vassiliev filtration has an alternative, more algebraic definition for many flavours of knot theory, for example braids and tangles, but notably not for knots. This view gives rise to connections between knot theory and quantum algebra. Finally, we review results -- many of them recent -- on extensions of non-numerical knot invariants to singular knots.
\end{abstract}

\maketitle

\section{Introduction}
\begin{figure}
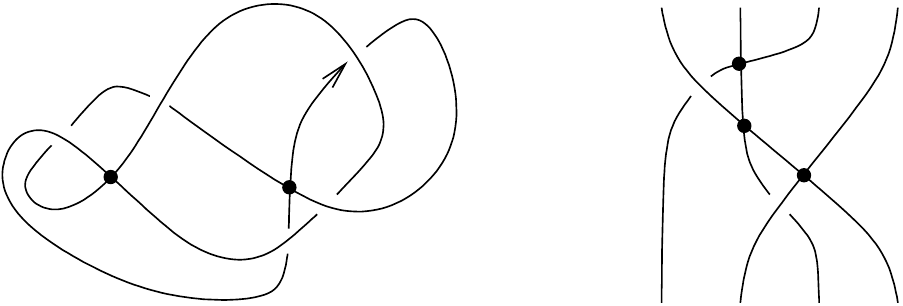
\caption{An example of a 2-singular knot -- that is, a knot with two double points; and a 3-singular braid on four strands.}\label{fig:SingKnotExample}
\end{figure}

\begin{figure}
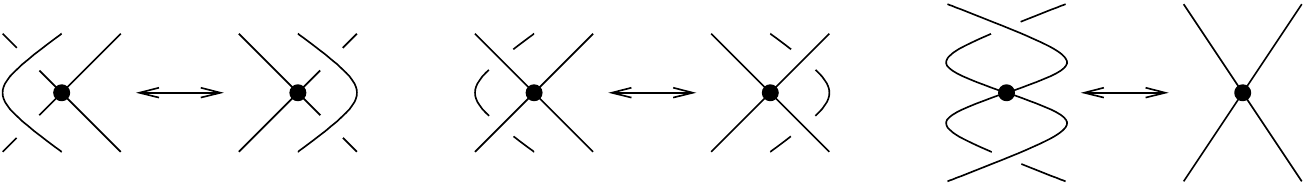
\caption{Reidemeister moves for singular knots, in addition to the usual Reidemeister 1, 2 and 3 moves.}
\label{fig:RMoves}
\end{figure}

A \textbf{singular knot} is an immersed circle in $\R^{3}$, whose singularities are limited to finitely many {\em transverse double points}. We will restrict our attention to {\em oriented} singular knots, which are equipped with a direction: Figure~\ref{fig:SingKnotExample} shows an example. Other knotted objects, such as braids, tangles, or knotted graphs have similar singular versions, see Figure~\ref{fig:SingKnotExample}
for an example of a singular braid. A \textbf{singular knot diagram} is a planar projection of a singular knot, which has two types of {\em crossings}: regular crossings with over/under strand information, and double points. Reidemeister moves can be formulated to describe singular knots modulo ambient isotopy, see Figure~\ref{fig:RMoves}. For minimal sets of oriented Reidemeister moves see \cite{BEHY}.

\section{The Vassiliev filtration}

Singular knots give rise to a decreasing filtration, called the {\em Vassiliev filtration} on the infinite dimensional vector space generated by isotopy classes of knots. Knot invariants which vanish on some step of this filtration are called {\em finite type} or {\em Vassiliev} invariants. Examples include many famous knot invariants, for example the coefficients of any knot polynomial -- after an appropriate variable substitution -- are of finite type \cite{BN}. In this section we describe the Vassiliev filtration for knots and other knotted objects. 

Let $\calK$ denote the set of (isotopy classes of) knots, and let $\Q\calK$ be the $\Q$-vector space\footnote{One may replace $\Q$ with one's favourite field of characteristic zero, at no cost.} spanned by elements of $\calK$. In other words, elements of $\Q\calK$ are formal linear combinations of knots. 

\begin{figure}
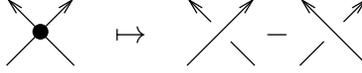
\caption{Resolution of singularities.}\label{fig:Resolve}
\end{figure}

Let $\calK^{n}$ denote the set of (isotopy classes of) $n$-singular knots, and  let $K\in \calK^{n}$. Each singularity of $K$ can be {\em resolved} two ways: by replacing the double point with an
over-crossing, or with and under-crossing. Note that the notion of ``over'' and ``under'' crossings don't depend on the choice of knot projection -- this is true not only in $\R^{3}$ but in any orientable manifold. We define a \textbf{resolution map} $\rho: \calK^{n} \to \Q\calK$ as follows: replace each singularity of $K$ by the difference of its two resolutions, as shown in Figure~\ref{fig:Resolve}. This produces a linear combination of $2^{n}$
knots with $\pm1$ coefficients, we call this the \textbf{resolution} of $K$. By an abuse of notation, we denote also by $\calK^{n}$ the linear span of the image $\rho(\calK^{n})$ in $\Q\calK$.

The subspaces $\calK^{n}$ are a decreasing filtration on $\Q\calK$, called the \textbf{Vassiliev filtration}:
\begin{equation}
\Q\calK=\calK^{0} \supset \calK^{1} \supset \calK^{2} \supset \calK^{3} \supset ...
\end{equation}
It's a worthwhile exercise to show that $\calK^{1}$ is the set of elements whose coefficients sum to zero: $$\calK^{1}=\{\sum_{i} \alpha_{i}K_{i} \in \Q\calK \, | \sum_{i} \alpha_{i}=0\}.$$

When one encounters a filtered space, a natural idea is to study its {\em associated graded} space instead: this enables inductive arguments and degree-by-degree computations. The associated graded space of knots equipped with the Vassiliev filtration is, by definition, $\mathcal A:= {\bigoplus}_{n=0}^{\infty} \calK^{n}/\calK^{n+1}$. The space $\mathcal A$ has a useful combinatorial--diagrammatic description in terms of {\em chord diagrams}. The study of finite type invariants is closely tied to the study of the space $\mathcal A$.

\begin{figure}
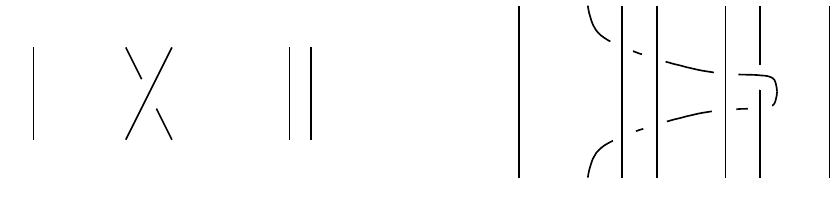
\caption{The generators $\sigma_{i}$ of $B_{n}$, and the generators $\sigma_{ij}$ of $PB_{n}$.}
\label{fig:BraidGens}
\end{figure}

We now examine the Vassiliev filtration from a more algebraic perspective through the example of braid groups. Recall that the \textbf{braid group} $B_{n}$ consists of
braids on $n$ strands up to braid isotopy, and the group operation is given by vertical stacking. The {\em Artin presentation} is a finite presentation for the group $B_{n}$ in terms of generators and relations (see also Figure~\ref{fig:BraidGens}):
$$B_{n}=\langle \sigma_{1}, ... , \sigma_{n-1} \, | \, \sigma_{i}\sigma_{i+1}\sigma_{i}=\sigma_{i+1}\sigma_{i}\sigma_{i+1}, \text{and} \, \,
\sigma_{i}\sigma_{j}=\sigma_{j}\sigma_{i} \, \, \text{when} \, |i-j|\geq 2 \rangle.$$

Following the strands of a braid from bottom to top determines a permutation of $n$ points; this gives rise to a group homomorphism $S: B_{n} \to S_{n}$ to the symmetric group. For example, (any resolution of) the braid of Figure~\ref{fig:SingKnotExample} is mapped to the cycle $(1324)$. The image $S(b)\in S_{n}$ of a braid $b$ is called the \textbf{skeleton} of $b$;
the kernel of this homomorphism is the \textbf{pure braid group} $PB_{n}$. The pure braid group also has a finite presentation due to Artin, in terms of the generators $\sigma_{ij} \, (i<j)$, shown in Figure~\ref{fig:BraidGens}, and somewhat more complicated relations. For more on pure braid presentations see \cite{MM}.

The Vassiliev filtration can be defined for braid groups the same way as it is defined for knots; we first discuss it for the pure braid group $PB_{n}$.
In this case $\Q PB_{n}$ is the {\em group algebra} of the pure braid group over the field $\Q$. Denote the linear subspace generated by the $\rho$-image of $n$-singular pure braids by $\calB^{n}$. Just like in the case of knots, 
\linebreak 
$\calB^{1}=\{\sum_{i} \alpha_{i}b_{i} \in \Q PB_{n} \, | \sum_{i} \alpha_{i}=0\}.$
This is a two-sided ideal in $\Q PB_{n}$, called the \textbf{augmentation ideal}.

However, in this case much more is true. Any $n$-singular braid can be written as a product of $n$ $1$-singular braids: one can ``comb'' the braid into $n$ horizontal levels so that each level has exactly one singularity. As a result, $\calB^{n}$ is simply the $n$-th power of the ideal $\calB^{1}$.

A similar statement holds for the braid group $B_{n}$, with the difference that in the definition of $\Q B_{n}$ linear combinations are only allowed within each coset of $PB_{n}$; in other words, only braids of the same skeleton can be combined. Steps of the Vassiliev filtration once again correspond to powers of the augmentation ideal. In fact, this phenomenon is much more general. Many flavours of knotted objects can be finitely presented as an algebraic structure of some kind: {\em tangles} form a {\em planar algebra}; {\em knotted trivalent graphs} have their own special set of operations \cite{Th}; virtual and
welded tangles form {\em circuit algebras} \cite{BD}. The Vassiliev filtration -- or its appropriate generalizations -- coincide with
the powers of the augmentation ideal in all of these examples.

\begin{figure}
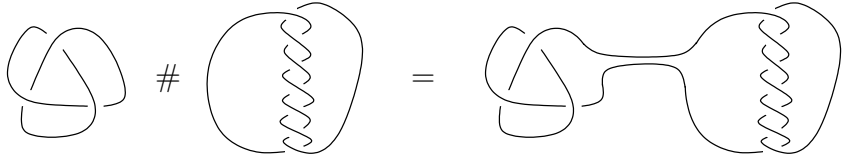
\caption{The connected sum operation of knots: the reader is encouraged to check that it is well-defined.}\label{fig:ConnSum}
\end{figure}

A notable exception is the case of knots. A ``multiplication operation'' does exist for knots: it is called \textbf{connected sum} and denoted $\#$; see Figure~\ref{fig:ConnSum} for a pictorial definition. It is not true, however, that any $n$-singular knot is the connected sum of $n$ 1-singular knots. Knots which cannot be expressed as a non-trivial connected sum are called \textbf{prime knots}, and there exist prime knots of arbitrary crossing number, in particular all {\em torus knots} are prime. In other words, knots are not {\em finitely generated} as an ``algebraic structure''.

The algebraic view outlined above gives rise to a deep connection between knot theory and quantum algebra. Let $\calK$ denote some class of knotted objects (such as knots, braids, tangles, etc). A central question in the study of finite type invariants is to find a {\em universal finite type invariant}, which contains all of the information that any finite type invariant can retain about  $\mathcal K$. More precisely, a universal finite type invariant is a filtered map $Z: \Q \mathcal K \to \hat{\mathcal A}$
which takes values in the {\em degree completed} associated graded space, and which satisfies a certain universality property. Perhaps the most famous universal finite type invariant is the {\em Kontsevich integral} of knots \cite{Ko}. It is still an open problem whether the Kontsevich integral separates knots.

Assume that some class of knotted objects $\mathcal K=\langle g_{1}, ..., g_{k}| R_{1},...,R_{l}\rangle$ forms a finitely presented algebraic structure with generators $g_{i}$ and relations $R_{j}$ (e.\ g.\ $\mathcal K$ may be the braid group). Assume one looks for a universal finite type invariant $Z$ which respects operations in the appropriate sense (e.\ g.\ $Z$ may be an algebra homomorphism). Then it is enough to find the values of the generators $Z(g_{i})\in \mathcal A$,
subject to equations arising from each $R_{j}$. In other words, one needs to solve a set of equations in a graded space. This set of equations often turns out to be interesting in its own right: for knotted trivalent graphs or {\em parenthesized tangles} they are the equations which define {\em Drinfeld associators} in quantum algebra \cite{MO, BN2, Da}; for {\em welded foams} they are the {\em Kashiwara-Vergne} equations of Lie theory \cite{BD}.

\section{Invariants of singular knots}

Vassiliev's idea in the early '90s was to extend number-valued knot invariants to singular knots using the resolution of singularities discussed above. This led to an explosion of activity in knot theory at the time. Since then singular knots have been studied in their own right by many researchers, and various non-numerical invariants have been extended to singular knots. Here we summarise some of these results.

\begin{figure}
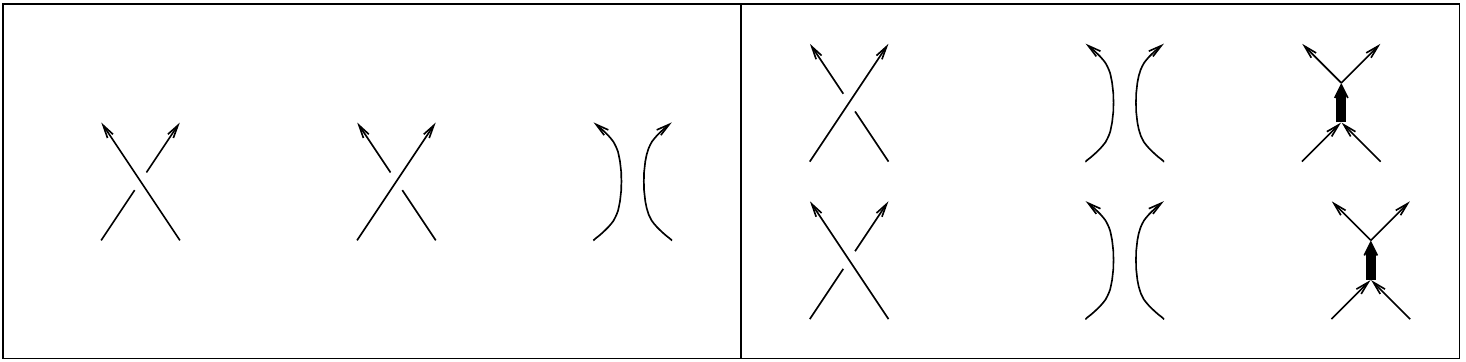
\caption{The original HOMFLY skein relation on the left; the HOMFLY skein relations for abstract singular knots \cite{MOY} on the right.}
\label{fig:HOMFLY}
\end{figure}

In \cite{MOY}, Murakami, Ohtsuli and Yamada developed a skein theory for the HOMFLY polynomial using a generalisation of singular knots called {\em abstract singular knots}. We present a simplified version of their result, which is used in several applications. The HOMFLY polynomial is determined by the skein relation shown on the left of Figure~\ref{fig:HOMFLY}, and its value on the unknot. 

\begin{figure}
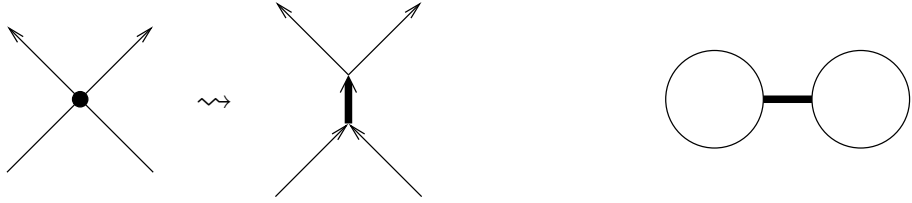
\caption{Creating an oriented abstract singular knot from a singular knot, and an abstract singular knot which cannot be oriented.}
\label{fig:ThickEdge}
\end{figure}

To define abstract singular knots, one replaces double points with an oriented {\em thick edge} from the incoming to the outgoing edges of the double point, creating an oriented trivalent graph embedded in $\R^{3}$, as in Figure \ref{fig:ThickEdge}. These trivalent graphs are called {\bf oriented abstract singular knots} and they are characterised by two properties:
\begin{enumerate}
\item Each vertex is incident to one thick and two thin edges.
\item At each vertex the thin edges are oriented the same way, while the thick edge is oriented oppositely.
\end{enumerate}

Embedded trivalent graphs which satisfy condition (1) are called {\bf abstract singular knots}. Not every such abstract singular knot can be oriented to satisfy (2), and hence they don't all arise from ordinary singular knots with double points: an example is shown in Figure~\ref{fig:ThickEdge}. 

Let $P_{n}$ be the one-variable spacialization of the HOMFLY polynomial with the substitution $x=q^{n}$, $y=q-q^{-1}$. In particular, $P_{0}$ is the Alexander polynomial, and $P_{2}$ is the Jones polynomial. Then $P_{n}$ satisfies the skein relations shown on the right of Figure~\ref{fig:HOMFLY}, reducing to the case of planar trivalent graphs, that is, abstract singular knots with no crossings. For these graphs \cite{MOY} provide a further set of skein relations which uniquely determine the value of $P_{n}$. 

Khovanov and Rozansky \cite{KR} use the \cite{MOY} calculus, and the homological algebra of {\em matrix factorizations}, in their categorification of the polynomials $P_{n}$. In \cite{OSS}, Ozsvath, Stipsitz and Szabo generalize knot Floer homology to abstract singular knots.

In \cite{KV} Kauffman and Vogel extend the Kauffman polynomial to singular knots -- in fact to knotted four-valent graphs with rigid vertices.
In \cite{JL} Juyumaya and Lambropoulou introduced a Jones-type invariant for singular knots, using the theory of singular braids and a Markov trace on Yokonuma--Hecke algebras. In \cite{Fi}, Fiedler extended the Kauffman state models of the Jones and Alexander polynomials to the context of singular knots.

Quandle-type invariants have been generalised and studied for singular knots and other types of knot-like objects (virtual knots, flat knots, pseudoknots). Authors who have contributed to this research include Churchill, Elhamdadi, Hajij, Henrich, Nelson, Oyamaguchi, and Sazdanovich \cite{HN,NOS,CEHN}. 

This short summary can not aim to be a comprehensive treatment of the rich body of research, spanning nearly three decades, on singular knots and related knotted objects. We hope that the reader will be inspired to explore some of the many pointers and references, and contribute to the future of the subject.

\end{document}